\documentclass[10pt]{article}

\newcommand{\given}{\big\vert}

\newcommand{\restrict}{\downharpoonright}

\newcommand{\WW}{\mathbb{W}}
\newcommand{\TT}{\mathbb{T}}

\newcommand{\PP}{\mathbb{P}}

\newcommand{\EE}{\mathbb{E}\,}

\newcommand{\cT}{\mathcal{T}}

\newcommand{\cF}{\mathcal{F}}
\newcommand{\cG}{\mathcal{G}}

\renewcommand{\Box}{\square}

\newcommand{\by}{\mathbf{y}}

\newcommand{\bP}{{\mathbf{P}}}

\newcommand{\bone}{\boldsymbol{1}}

\newcommand{\tTT}{{\tilde{\mathbb{T}}}}

\newcommand{\tmu}{{\tilde{\mu}}}

\usepackage{amsfonts}
\usepackage{amsmath}
\usepackage{a4}
\usepackage{amssymb}
\usepackage{theorem}

\newtheorem{theorem}{Theorem}
\newtheorem{lemma}[theorem]{Lemma}
\newtheorem{corollary}[theorem]{Corollary}
{\theorembodyfont{\rmfamily}}
{\theorembodyfont{\rmfamily}}

\begin{document}

\author{James B.\ Martin}
\title{Reconstruction Thresholds on Regular Trees}
\date{CNRS and Universit\'e Paris 7}
\maketitle

%

\begin{abstract}
We consider a branching random walk with binary state space 
and index set $\TT^k$, the infinite rooted tree in which each node
has $k$ children
(also known as the model of \textit{broadcasting on a tree}).
The root of the tree takes 
a random value $0$ or $1$,
and then each node passes a value
independently to each of its children 
according to a $2\times2$ transition matrix $\bP$.
We say that \textit{reconstruction is possible}
if the values at the $d$th level of the tree
contain non-vanishing information about the value
at the root
as $d\to\infty$.
Adapting a method of Brightwell and Winkler,
we obtain new conditions under which reconstruction
is impossible, both in the general case
and in the special case $p_{11}=0$.
The latter case is closely related to the 
\textit{hard-core model} from statistical physics;
a corollary of our results is that, for 
the hard-core model on the $(k+1)$-regular tree
with activity $\lambda=1$,
the unique simple invariant Gibbs measure
is extremal in the set of Gibbs measures, for any $k$.
\end{abstract}

\tableofcontents

\section{Introduction}
\label{introduction}

\subsection{Branching random walk model}
\label{BRW}

We consider a model of a branching random walk (BRW) 
indexed by the rooted tree
$\TT^k$, in which every node has
$k$ children. 

Let $\bP=\{p_{ij}, i,j=0,1\}$ 
be a $2\times2$ stochastic matrix, 
which we regard as a transition matrix on the set $\{0,1\}$.
Each node $u\in\TT^k$ will carry 
a value $\phi(u)\in\{0,1\}$, generated as follows.
The root takes value 0 with
probability $\pi_0=p_{01}/(p_{01}+p_{10})$
and value 1 with probability $\pi_1=1-\pi_0$.
Thereafter the configuration on $\TT^k$ 
is generated recursively;
if a node has value $i\in\{0,1\}$,
each of its $k$ children takes
the value 0 with probability $p_{i0}$ 
and the value 1 with probability $p_{i1}$,
all choices being made independently.

We write $\phi=\{\phi(u), u\in\TT^k\}$
for a configuration on the whole tree, and denote by
$\mu$ the probability measure on $\{0,1\}^{\TT^k}$
resulting from this branching random walk construction.

For a node $u\in\TT^k$,
let $\TT^k(u)$ be the subtree consisting of $u$
and all its descendants. By the choice of $\pi_0$,
we have a translation invariance property for $\mu$;
namely that $\mu(\phi(u)=0)=\pi_0$ for every $u\in\TT^k$,
and so 
for any $u$, $v\in\TT^k$, the configurations on
$\TT^k(u)$ and $\TT^k(v)$ have the
same distribution,
under a natural mapping between the subtrees $\TT^k(u)$
and $\TT^k(v)$.

We are interested in the following question of 
\textit{reconstruction}: for $d\geq1$,
how much information about the value at node $u$
is given by the values of the $d$th generation
of its descendants?

Questions of this sort arise in several 
contexts -- for example genetics, communication theory 
and statistical physics -- and have
been quite widely studied in the last few years; 
see Mossel \cite{Mosselsurvey} for a survey,
and \cite{EKPS, Bleher, Ioffe, KMP, Mossel2nd, MosPer, BriWinhard, JanMos}
for a variety of approaches to this sort of model
(which can of course be considerably generalised
from our particular setting of a binary state space 
and a regular tree).

The question above can be made precise in several
(often equivalent) ways. We use the following formulation.

Let $\WW_d(u)$ be the set of descendants of $u$ at 
distance exactly $d$ from $u$.
For a set $S\subseteq\TT^k$, write $\sigma(S)$
for the $\sigma$-algebra of events which depend only
on the values $\{\phi(u),u\in S\}$.

Define the random variable 
\[
A(d,u)=\mu\big(\phi(u)=0 \given \sigma(\WW_d(u))\big),
\]
that is, the conditional probability that the value at $u$ is 
0, given only the information from the $d$th generation of its 
descendants.

From the 
indepdence structure given by the 
branching random walk construction,
additional
knowledge of any information from nodes beyond the $d$th generation
does not change the conditional distribution of the value of $u$;
that is,
\[
A(d,u)=\mu\Big(\phi(u)=0 \given 
\sigma\big(\bigcup_{d'=d}^{\infty}\WW_{d'}(u)\big)\Big).
\]

Of course, if $d_1>d_2$, then
\[
\sigma\big(\bigcup_{d'=d_1}^{\infty}\WW_{d'}(u)\big)
\subseteq
\sigma\big(\bigcup_{d'=d_2}^{\infty}\WW_{d'}(u)\big),
\]
so by the backwards martingale convergence theorem
(see e.g.\ Section 14.4 of \cite{Williams}), 
we have that $A(d,u)\to A(u)$ a.s.\ 
as $d\to\infty$, where
\[
A(u)=\mu(\phi(u)=0 \given \cT(u));
\]
here $\cT(u)$ is the \textit{tail $\sigma$-algebra} of
descendants of $u$,
defined by 
\[
\cT(u)=\bigcap_{d=1}^{\infty}
\sigma\big(\bigcup_{d'=d}^{\infty}\WW_{d'}(u)\big).
\]

By the translation invariance property above,
the random variable
$A(u)$ has the same distribution for all $u\in\TT^k$.

\medskip

\noindent\textit{Definition:} 
We say that \textit{reconstruction is impossible}
(for a given $\bP$ and $k$)
if the random variable $A(u)$ is 
almost surely constant, and otherwise that 
\textit{reconstruction is possible}.

\medskip

A complete answer to the question of when reconstruction is possible
is currently only known for the case where $\bP$ is symmetric.
Then let $p_{00}=p_{11}=1-\epsilon$; reconstruction is 
possible if and only if $k(1-2\epsilon)^2>1$
(see for example \cite{Bleher, EKPS, Ioffe}).

In general, however, there are gaps between the best known 
necessary and sufficient conditions for reconstruction to be possible.
In this paper we give new conditions on $\bP$ under which we show
that reconstruction is impossible. 

In Proposition 4.1 of \cite{MosPer},
Mossel and Peres show that
reconstruction is impossible whenever 
\begin{equation}
\label{mospercond}
\frac{(p_{10}-p_{00})(p_{01}-p_{11})}{\min\{p_{00}+p_{10},p_{01}+p_{11}\}}
\leq \frac{1}{k}.
\end{equation}
(The two factors in the numerator are the same; we write
this form simply to preserve the symmetry between the states 0 and 1).
We improve the bound to give the following condition:

\begin{theorem}
\label{generaltheorem}
Reconstruction is impossible whenever
\begin{equation}
\label{generalcond}
\Big(\sqrt{p_{00}p_{11}}-\sqrt{p_{01}p_{10}}\Big)^2\leq\frac{1}{k}.
\end{equation}
\end{theorem}

A calculation shows that the LHS of (\ref{generalcond})
is always less than or equal to that of (\ref{mospercond}),
with equality in the following special cases:
(i) $\bP$ is symmetric; (ii) $p_{ij}=0$ for some $i,j$;
(iii) $p_{00}=p_{10}$, $p_{01}=p_{11}$. Note that for
symmetric $\bP$, (\ref{generalcond}) becomes
the condition that $k(1-2\epsilon)^2\leq 1$,
and our proof of Theorem \ref{generaltheorem} 
gives another proof that reconstruction is impossible
under this condition.

We then focus on the special case where $p_{11}=0$ (of course,
the case $p_{00}=0$ is analogous). This case is closely
related to the \textit{hard-core model} from statistical physics,
and has been recently studied by 
Brightwell and Winkler \cite{BriWinhard} 
and Rozikov and Suhov \cite{RozSuh}.
Certain specific properties in this case allow a
more sophisticated argument which gives 
a much better condition than is obtained by putting $p_{11}=0$ in 
Theorem \ref{generaltheorem}.

\subsection{Hard-core model}
\label{hardcore}
In this section we state our result for the case $p_{11}=0$
and explain the correspondence with 
the hard-core model on a regular tree.

Following \cite{BriWinhard},
we parametrise $\bP$
by the quantity $w>0$, setting
\begin{equation}
\label{matrix}
\bP= \begin{pmatrix} p_{00} & p_{01} \\ p_{10} & p_{11} \end{pmatrix}
= \begin{pmatrix} \frac{1}{1+w} & \frac{w}{1+w} \\ 1 & 0 \end{pmatrix},
\end{equation}
or equivalently by the quantity $\lambda=w(1+w)^k>0$,
whose significance we explain later;
note that the correspondence between $\lambda>0$ and $w>0$ is
one-to-one and monotonic.

Let $\lambda_c=\lambda_c(k)$ be the infimum of the set of $\lambda$
such that reconstruction is possible. If follows 
from Proposition 12 of \cite{Mossel2nd} that in fact
reconstruction is possible for any $\lambda>\lambda_c$
(so that $\lambda_c$ is also the supremum of the set of $\lambda$
such that reconstruction is impossible).

Brightwell and Winkler \cite{BriWinhard}
show that, as $k\to\infty$,
\begin{equation}\label{BriWin}
\frac{1+o(1)}{\ln k}\leq \lambda_c(k) \leq (\ln k)^2\big(1+o(1)\big).
\end{equation}
We improve the lower bound to give the following:
\begin{theorem}
\label{hardtheorem}
$\lambda_c(k)>e-1$ for all $k$.
\end{theorem}

(For the equivalent threshold value $w_c$ 
with $w_c(1-w_c)^k=\lambda_c$, one can deduce that
$w_c(k)>(\ln k-\ln\ln k)/k$ for all $k$).

We will now describe the correspondence 
between the BRW and the hard-core model, 
and explain (without proofs) the significance of Theorem \ref{hardtheorem}
for the hard-core model on the $(k+1)$-regular tree.
For more details on the correspondence between the two models, 
see also \cite{BriWinhard} and its references.

We denote the $(k+1)$-regular tree
by $\tTT^k$. 
We can still regard
$\tTT^k$ as a rooted tree, in which the root node has $k+1$
children and every other node has $k$ children.
We can then carry out the branching random walk construction
on $\tTT^k$ in exactly the same way as we did on $\TT^k$; now the root
has $k+1$ rather than $k$ children, but the values at these $k+1$ children
are chosen i.i.d.\ according to the value at the root and 
the transition matrix $\bP$ just as before. We will write $\tmu$ for the 
probability measure on $\{0,1\}^{\tTT^k}$ resulting from this construction.

The independence structure of the random walk implies that the measure 
$\tmu$ is 
\textit{simple}, by which we mean that, for any $u$,
the configurations
\[
\{\phi(v), v\in C_1(u)\}, \ldots, \{\phi(v), v\in C_{k+1}(u)\}
\]
are mutually independent given $\phi(u)$,
where the $C_i$ are the connected components of $\tTT^k\setminus\{u\}$. 
Although we have defined $\tmu$ in an asymmetric way, it's also 
the case that 
it is 
\textit{invariant},
in the sense that it is preserved by any automorphism of $\tTT^k$.
In particular, the choice of the root is
not important.

To introduce the hard-core model, we first consider
the case of a finite graph with node-set $S$ (and 
some neighbour relation).

We can identify a configuration $\phi\in\{0,1\}^S$
with the subset 
$I_\phi:=$$\{u\in\tTT^k:$$\phi(u)=1\}$ of $S$.

A set $I\subseteq S$ is called an \textit{independent set}
if no two neighbours in the graph are both members of $I$.

The \textit{hard-core measure} on $S$ with \textit{activity} $\lambda>0$
is the probability measure $\nu$ on $\{0,1\}^S$
such that 
\[
\nu\big(
I_\phi \text{ is an independent set }\big)=1,
\]
and such that for an independent set $I_0$,
$\nu(I_\phi=I_0)$ is proportional to $\lambda^{|I_0|}$.
Thus in fact
\[
\nu(\phi=\phi_0)=Z_\lambda^{-1}\lambda^{|I_{\phi_0}|}
1\big(I_{\phi_0} \text{ is an independent set }\big),
\]
where we have the normalising factor
\[
Z_\lambda=\sum_{\phi_0:I_{\phi_0}\text{ is independent}}\lambda^{|I_{\phi_0}|}.
\]
When $\lambda=1$, $I_\phi$ has the uniform distribution
over the set of independent subsets of $S$.

An equivalent characterisation is that 
$\nu$ is the unique probability measure such that,
for any $\phi_0\in\{0,1\}^S$ and any $u\in S$,
\begin{equation}
\label{charact}
\nu\Big(
\phi(u)=1 \given \phi(v)=\phi_0(v) \text{ for all } v\neq u
\Big)
=
\frac{\lambda}{1+\lambda}
\bone \left\{
I_{\phi_0}\cup\{u\} \text{ is independent}
\right\}
\end{equation}

The condition (\ref{charact}) 
makes sense equally when $S$ is infinite, 
except that (since conditional probabilities 
are only well defined up to almost sure equality)
we should now only demand the condition holds 
for $\nu$-almost all $\phi_0$. Putting $S=\tTT^k$,
we say that a probability measure $\nu$ satisfying 
(\ref{charact}) 
(for all $u\in\tTT^k$ and $\nu$-almost all $\phi_0$)
is a \textit{Gibbs measure}
for the hard-core model on $\tTT^k$ with activity $\lambda$.

It is quite straightforward, 
to show that the measure $\tmu$ defined above by 
the BRW construction with $\bP$ as in
(\ref{matrix}) is a Gibbs measure for 
the hard-core model with activitiy $\lambda$.
However, now that the state space is infinite,
it's no longer the case that such a measure need be unique.
In fact, there is a critical point 
$\lambda'_c=$$\lambda'_c(k)=$$k^k/(k-1)^{(k+1)}$
(identified by Kelly \cite{Kellyhard});
for $\lambda\leq \lambda'_c$,
$\tmu$ is the \textit{only} Gibbs measure,
whereas for $\lambda>\lambda'_c$, there are others.

Nevertheless, for any $\lambda$, 
$\tmu$ is the only \textit{simple invariant} Gibbs measure.

The set of Gibbs measures forms
a simplex; 
that is, any mixture of Gibbs measures is also a Gibbs measure,
and in particular 
there is a set of \textit{extremal}
Gibbs measures such that every Gibbs measure
is expressible in a unique way as a mixture of extremal measures.
For $\lambda > \lambda'_c$,
we can therefore ask 
whether the measure $\tmu$ is extremal
(equivalently, not expressible as a mixture of other Gibbs measures).

It turns out that $\tmu$ is extremal at activity $\lambda$
if and only if reconstruction is impossible for 
the corresponding branching random walk on $\TT^k$ with
transition matrix $\bP$. 
(This is a consequence of the general fact that a Gibbs
measure is extremal iff it is trivial on the tail $\sigma$-algebra,
and of the independence structure
given by the BRW constructions of $\mu$ and $\tmu$).
Hence the reconstruction threshold $\lambda_c$ 
defined after (\ref{matrix})
is also the \textit{extremality} threshold for $\tmu$;
Theorem \ref{hardtheorem} 
therefore shows that 
whenever $\lambda\leq e-1$, the unique simple invariant
Gibbs measure $\tmu$ for the hard-core model with activity $\lambda$
is extreme, for any $k$.

In particular, $\tmu$ is extreme in the special case $\lambda=1$
for any $k$. 

\subsection{Outline of proof}
\label{outline}
Our approach to proving Theorems \ref{generaltheorem} and \ref{hardtheorem} 
is closely related to the method used by Brightwell
and Winkler to prove the lower bound in (\ref{BriWin})
for the hard-core model
\cite{BriWinhard}. 

We first develop a coupling between the
distributions of the random variable $A(u)$ 
conditioned on two different events,
with certain additional properties beyond those 
used in \cite{BriWinhard}. We then use 
this coupling to establish 
a recursion linking the distribution of $A(u)$ 
to those of $A(u_1),\ldots,A(u_k)$, where $u_1,\ldots, u_k$
are the children of $u$. 
(Of course, we already know from the translation 
invariance property described in Section \ref{BRW} 
that all of these distributions are the same).
If the recursion relation 
is contractive in a suitable sense, we obtain
that $A(u)$ must be a.s.\ constant.

In Section \ref{conditional},
we first prove a lemma on conditional probabilities
in a more general setting. 
Specialising to our context, we obtain the existence
of a coupling of a pair of random variables $A_0$, $A_1$
with the following properties:
\begin{itemize}
\item[(i)]
The distribution of $A_0$ 
is the distribution of $A(u)$ under $\mu$ 
conditioned on the event $\{\phi(u)=0\}$;
\item[(ii)]
The distribution of $A_1$ 
is the distribution of $A(u)$ under $\mu$ 
conditioned on the event $\{\phi(u)=1\}$;
\item[(iii)] 
With probability 1,
either $A_0=A_1$ or $A_1\leq\pi_0\leq A_0$;
\item[(iv)]
If $A_0=A_1$ with probability 1,
then both are equal to $\pi_0$ with probability 1,
and so also $A(u)=\pi_0$ a.s.\ under $\mu$.
\end{itemize}

We develop the recursion relations and 
complete the proofs in Section \ref{mainproof}.

The full properties of the coupling are only needed
in the hard-core case, where a
particular convexity property holds 
for the recursion relations.
The argument in the general case is not as powerful, 
and rather than all of property (iii) above, we use only
that $A_1\leq A_0$ with probability 1. 
Restricting the bound in Theorem \ref{generaltheorem}
to the case $p_{11}=0$ gives 
a much weaker bound than that in Theorem \ref{hardtheorem}
(in fact, one obtains only the bound $\lambda_c\geq\lambda'_c$
where $\lambda'_c$ is the threshold for the uniqueness
of the Gibbs measure; this bound is obvious 
in the context of the hard-core model since
if the Gibbs measure is unique it is trivially extreme).

\section{Conditioned conditional probabilities}
\label{conditional}
We first consider the setting of a general
probability space
$(\Omega, \cF, \PP)$.
Let $B\in\cF$ be an event with probability $\pi_0$=$1-\pi_1$,
and suppose $0<\pi_0<1$. Write $B^C$ for the complement of $B$.
Let $\cG$ be a sub-$\sigma$-algebra of $\cF$.

We consider the random variable $\PP(B \given \cG)$
(which is the $\cG$-measurable 
random variable, unique up to almost sure equality,
such that for all $D\in\cG$
\begin{equation}\label{cpdef}
\PP(D\cap B)=\int_D \PP(B \given \cG)(\omega)d\PP(\omega).
\end{equation}
See for example Chapter 9 of \cite{Williams}
for background on conditional probabilities).

\begin{lemma}
Suppose $0\leq p_0\leq p_1\leq 1$,
and that $D\in\cG$ with
\begin{equation}\label{Dcond}
\PP(B \given \cG)(\omega)\in[p_0,p_1]
\text{ for all } \omega\in D.
\end{equation}
Then
\[
\frac{\pi_1}{\pi_0}\frac{p_0}{1-p_0}\PP(D \given B^C)
\leq \PP(D \given B)
\leq
\frac{\pi_1}{\pi_0}\frac{p_1}{1-p_1}\PP(D \given B^C).
\]
\end{lemma}

\noindent\textit{Proof:}
From (\ref{cpdef}) and (\ref{Dcond}) we have
\begin{gather*}
p_0\PP(D)\leq\PP(D\cap B)\leq p_1\PP(D)
\\
\intertext{and}
(1-p_1)\PP(D)\leq\PP(D\cap B^C)\leq (1-p_0)\PP(D).
\end{gather*}
Combining these we get
\[
\frac{p_0}{1-p_0}\PP(D\cap B^C)
\leq \PP(D\cap B)
\leq \frac{p_1}{1-p_1}\PP(D\cap B^C).
\]
Since $\PP(D \given B)=\PP(D\cap B)/\pi_0$ 
and $\PP(D \given B^C)=\PP(D\cap B^C)/\pi_1$,
the result follows.
$\hfill\Box$

\bigskip

In particular, if $J$ is a subset of the interval $[0,\pi_0)$
then we can set
$D=\{\omega:\PP(B\given\cG)(w)\in J\}$
to obtain
\[
\PP\Big\{\PP(B \given \cG)\in J \given B\Big\}
\leq
\PP\Big\{\PP(B \given \cG)\in J \given B^C\Big\},
\]
while if $J\subseteq(\pi,1]$ then the inequality is
reversed. In each case equality holds only if both sides are 0.
Hence:

\begin{corollary}
There exists a coupling of two random variables $Y_0$ 
and $Y_1$, such that $Y_0$ has the distribution of $\PP(B \given \cG)$
conditioned on $B$ occurring, such that $Y_1$ has the distribution
of $\PP(B \given \cG)$ conditioned on $B$ not occurring, and such that:
\begin{itemize}
\item[(i)]whenever $Y_0<\pi_0$, then $Y_1=Y_0$, and
\item[(ii)]whenever $Y_1>\pi_0$, then $Y_1=Y_0$.
\end{itemize}
Therefore either $Y_0=Y_1$ or $Y_1\leq \pi_0 \leq Y_0$.

Also the distributions of $Y_0$ and $Y_1$ are identical iff
$Y_0=Y_1=\pi_0$ with probability 1, or equivalently 
iff $\PP(B \given \cG)=\pi_0$ with probability 1.
\end{corollary}

Applying this result with $B=\{\phi(u)=0\}$, with $\cG=\cT(u)$,
with $\PP=\mu$  
and so with $A(u)=\PP(B \given \cG)$, we obtain the coupling 
of $A_0, A_1$ with the properties
claimed in 
Section \ref{outline}.
 
\section{Recurrences for likelihood ratios}
\label{mainproof}

Let $u\in\TT^k$ and 
let $\by$ be a configuration on the set $\WW_d(u)$ 
(the descendants of $u$ at distance exactly $d$).

For $S\subset \TT^k$,
write $\phi\restrict_S$ for the configuration 
$\phi$ restricted to $S$.

Define the ``likelihood functions''
\begin{align*}
q^{(0)}(d,\by)=\mu\big(\phi\restrict_{\WW_d(u)}=\by  \given  \phi(u)=0)
\\
q^{(1)}(d,\by)=\mu\big(\phi\restrict_{\WW_d(u)}=\by  \given  \phi(u)=1).
\end{align*}
For $i=0,1$, the function $q^{(i)}(d,\by)$ gives the 
probability of observing the configuration $\by$ on 
the set of descendants of $u$ at distance $d$,
given that the value at $u$ itself is $i$.
(Note that because of the translation invariance property
noted in Section \ref{BRW}, the choice of $u$ is not important).

Define also the ``likelihood ratio'' function
\[
q(d,\by)=\frac{q^{(0)}(d,\by)}{q^{(1)}(d,\by)}.
\]

Let $d\geq2$ and
let the children of $u$ be $u_1,\ldots,u_k$.
A configuration $\by$ on $\WW_d(u)$ corresponds to a set of
configurations $\by_1,\ldots,\by_k$
on $\WW_{d-1}(u_1),\ldots,\WW_{d-1}(u_k)$.
We then have
\begin{gather}
\nonumber
q^{(0)}(d,\by)=\prod_{j=1}^k
\big[
p_{00}q^{(0)}(d-1,\by_j)
+p_{01}q^{(1)}(d-1,\by_j)
\big]
\\
\nonumber
q^{(1)}(d,\by)=\prod_{j=1}^k
\big[
p_{10}q^{(0)}(d-1,\by_j)
+p_{11}q^{(1)}(d-1,\by_j)
\big]
\\
\intertext{and so}
\label{qrecurrence}
\begin{split}
q(d,\by)
&=\prod_{j=1}^k
\left\{
\frac{
p_{00} q^{(0)}(d-1,\by_j) 
+
p_{01} q^{(1)}(d-1,\by_j) 
}
{
p_{10} q^{(0)}(d-1,\by_j) 
+
p_{11} q^{(0)}(d-1,\by_j) 
}
\right\}
\\
&=
\left(\frac{p_{00}}{p_{10}}\right)^k
\prod_{j=1}^k
\left\{
1+\frac{c_0-c_1}{q(d-1,\by_j)+c_1}
\right\},
\end{split}
\end{gather}
where we define 
\[
c_0=\frac{p_{01}}{p_{00}}, 
\quad 
c_1=\frac{p_{11}}{p_{10}}.
\]

Define also 
$a(d,\by)=\mu\left(\phi(u)=0 \given \phi\restrict_{\WW_d(u)}=\by\right)$.
The function $a$ gives the conditional probability that
the value at the node $u$ is 0, given a configuration on the 
set of descendants of $u$ at distance $d$.
We have
\begin{align}
\nonumber
a(d,\by)
&=\frac{\pi_0 q^{(0)}(d,\by)}{\pi_0 q^{(0)}(d,\by)+\pi_1 q^{(1)}(d,\by)}
\\
\label{aqrel}
&=\frac{1}{1+\frac{\pi_1}{\pi_0 q(d,\by)}}.
\end{align}

Returning to the random variable 
$A(d,u)=\mu\left(\phi(u)=0\given\sigma(\WW_d(u))\right)$
defined in Section \ref{introduction},
we have 
\[
A(d,u)=a\left(d,\phi\restrict_{\WW_d(u)}\right),
\]
that is, the function $a(d,.)$ applied to the 
actually observed values of the configuration $\phi$ on $\WW_d(u)$.

Similarly define
\[
Q(d,u)=q\left(d,\phi\restrict_{\WW_d(u)}\right).
\]
From (\ref{aqrel}),
we have
\[
A(d,u)=\frac{1}{1+\frac{\pi_1}{\pi_0 Q(d,u)}},
\;\;\;
Q(d,u)=\frac{\pi_1}{\pi_0}\left(\frac{1}{1-A(d,u)}-1\right).
\]
Recalling $A(d,u)\to A(u)$ a.s.,
we have $Q(d,u)\to Q(u)$ a.s., where
\[
Q(u)=\frac{\pi_1}{\pi_0}\left(\frac{1}{1-A(u)}-1\right).
\]
From (\ref{qrecurrence}) we get
\begin{gather*}
Q(d,u)=\left(\frac{p_{00}}{p_{10}}\right)^k
\prod_{j=1}^k\left(1+\frac{c_0-c_1}{Q(d-1,u_j)+c_1}\right),
\\
\intertext{and, taking $d\to\infty$,}
Q(u)=\left(\frac{p_{00}}{p_{10}}\right)^k
\prod_{j=1}^k\left(1+\frac{c_0-c_1}{Q(u_j)+c_1}\right).
\end{gather*}

Now put 
\begin{gather*}
L(u,d)=\ln Q(u,d)
\\
\intertext{and}
\begin{split} 
L(u)&=\ln Q(u)
\\
&=\ln\left[\frac{\pi_1}{\pi_0}\left(\frac{1}{1-A(u)}-1\right)\right].
\end{split}
\end{gather*}
We have $L(u,d)\to L(u)$ a.s.\ as $d\to\infty$,
and
\begin{equation}
\label{Lrecurrence}
L(u)=k\ln\left(\frac{p_{00}}{p_{10}}\right)
+\sum_{j=1}^k \ln\left(1+\frac{c_0-c_1}{\exp\big(L(u_j)\big)+c_1}\right).
\end{equation}

Since $L(u)$ can be written as a strictly increasing function
of $A(u)$, with $A(u)=\pi_0$ corresponding to $L(u)=0$,
we can translate the coupling of $A_0, A_1$ 
described in Section \ref{outline} and
proved in Section \ref{conditional} into a coupling
of two random variables $L_0, L_1$ 
with the following properties:
\begin{itemize}
\item[(i)]
The distribution of $L_0$ 
is the distribution of $L(u)$ conditioned on the event $\{\phi(u)=0\}$;
\item[(ii)]
The distribution of $L_1$ 
is the distribution of $L(u)$ conditioned on the event $\{\phi(u)=1\}$;
\item[(iii)] 
With probability 1,
either $L_0=L_1$ or $L_1\leq 0 \leq L_0$;
\item[(iv)]
If $L_0=L_1$ with probability 1,
then both are equal to $0$ with probability 1,
and then also $A(u)=\pi_0$ with probability 1.
\end{itemize}

So to conclude that $A(u)$ is a.s.\ constant,
it's enough to show that $\EE (L_0-L_1)=0$.

Note 
that from (\ref{Lrecurrence})
it's immediate that the distribution of $L(u)$ has compact support; 
hence the same is true for $L_0$ and $L_1$, 
and certainly 
$\EE |L_0|< \infty$, $\EE |L_1|<\infty$.

Again let $u_1,\ldots,u_k$ be the children of a node $u$.
From the branching random walk construction we get
the following information.

\noindent
Conditional on $\phi(u)=0$:
\begin{itemize}
\item[]
the $\phi(u_j), j=1,\ldots,k$
are i.i.d.\ taking value 0 w.p.\ $p_{00}$ 
and value 1 w.p.\ $p_{01}$. Then the 
$L(u_j)$ are i.i.d., and the distribution of each 
is a mixture of the distribution of $L_0$ 
(with weight $p_{00}$) and the 
distribution of $L_1$ (with weight $p_{01}$).
\end{itemize}

\medskip
\noindent
Conditional on $\phi(u)=1$:
\begin{itemize}
\item[]
the $\phi(u_j), j=1,\ldots,k$
are i.i.d.\ taking value 0 w.p.\ $p_{10}$ 
and value 1 w.p.\ $p_{11}$. Then the 
$L(u_j)$ are i.i.d., and the distribution of each 
is a mixture of the distribution of $L_0$ 
(with weight $p_{10}$) and the 
distribution of $L_1$ (with weight $p_{11}$).
\end{itemize}

Hence from (\ref{Lrecurrence}),
\begin{gather}
\nonumber
\EE L_0=
-k\ln\left(\frac{p_{00}}{p_{10}}\right)
+k\left[
p_{00}\EE\ln\left(1+\frac{c_0-c_1}{\exp(L_0)+c_1}\right)
+p_{01}\EE\ln\left(1+\frac{c_0-c_1}{\exp(L_1)+c_1}\right)
\right].
\\
\intertext{and}
\nonumber
\EE L_1=
-k\ln\left(\frac{p_{00}}{p_{10}}\right)
+k\left[
p_{10}\EE\ln\left(1+\frac{c_0-c_1}{\exp(L_0)+c_1}\right)
+p_{11}\EE\ln\left(1+\frac{c_0-c_1}{\exp(L_1)+c_1}\right)
\right].
\\
\intertext{Subtracting and using the fact that $p_{00}-p_{10}=p_{11}-p_{01}$,
we obtain}
\label{expect}
\EE (L_0-L_1)
=k\EE\left[
f(L_0)-f(L_1)
\right],
\\
\intertext{where}
\label{fdef}
f(x)=\left(p_{11}-p_{01}\right)
\ln\left(1+\frac{c_0-c_1}{e^x+c_1}\right).
\end{gather}

\subsection{General case}
\label{generalproof}
If $p_{01}=p_{11}$, then also $c_0=c_1$,
and the function $f$ is constant. In that case,
(\ref{expect}) shows that $L_0=L_1$ with probability 1,
and reconstruction is impossible.

So assume that $p_{01}\neq p_{11}$. 
Then also $c_0\neq c_1$, and
is easy to check that the function $f$ defined in (\ref{fdef}) is 
always strictly increasing
(since $c_0>c_1$ iff $p_{01}>p_{11}$).

A calculation gives that $\sup f'(x)$ is attained exactly when
\[
x=\frac{1}{2}\Big[
\ln \big(p_{01}p_{11}\big) -\ln \big(p_{00}p_{10}\big)
\Big],
\]
and that in fact
\begin{equation}\label{sup}
\sup f'(x)=\Big(\sqrt{p_{00}p_{11}}-\sqrt{p_{01}p_{10}}\Big)^2.
\end{equation}

Since 
we know that $L_1\leq L_0$ with probability 1,
we then have that 
\[
0\leq f(L_0)-f(L_1)\leq \sup_{x} f'(x) (L_0-L_1)
\]
with equality on the RHS iff $L_0=L_1$ with probability 1.
Hence, from (\ref{expect}), 
\[
\EE(L_0-L_1)\leq k\sup_x f'(x) \EE(L_0-L_1),
\]
with equality iff both sides are $0$.
So to show that $\EE(L_0-L_1)=0$, 
and therefore that reconstruction is impossible,
it's enough to show that 
\[
k\sup f'(x)\leq 1.
\]

Using (\ref{sup}), we see that
(\ref{generalcond}) indeed implies that 
reconstruction is impossible, and the proof 
of Theorem \ref{generaltheorem} is done.

\subsection{Hard-core case}
\label{hardproof}
Recall that 
\[
\begin{pmatrix} p_{00} & p_{01} \\ p_{10} & p_{11} \end{pmatrix}
= \begin{pmatrix} \frac{1}{1+w} & \frac{w}{1+w} \\ 1 & 0 \end{pmatrix};
\]
we have also that $\pi_0=(1+w)/(1+2w)$, $\pi_1=w/(1+2w)$,
and $c_0=w$, $c_1=0$. We have also defined $\lambda=w(1+w)^k$.

Equation (\ref{fdef}) now becomes
\[
f(x)=-\frac{w}{1+w}\ln\left(1+we^{-x}\right),
\]
and now the function $f$ is concave as well as 
strictly increasing.
Hence in particular, if
$x_0<x_1$ and $y_0<y_1$ with $x_0\leq y_0$
and $x_1\leq y_1$, then
\begin{equation}
\label{concave}
0
\leq
\frac{f(y_1)-f(y_0)}{y_1-y_0}
\leq
\frac{f(x_1)-f(x_0)}{x_1-x_0}.
\end{equation}

Recurrence (\ref{Lrecurrence})
now becomes
\[
L(u)=-k\ln(1+w)+\sum_{j=1}^k \ln\left(1+we^{-L\left(u_j\right)}\right),
\]
and so in particular $L(u)$ is always greater than or equal 
to $-k\ln(1+w)$; the same is therefore true of $L_0$ and $L_1$ also.

Combining this with property (iii) of the coupling
described after (\ref{Lrecurrence}),
we in fact have that 
with probability 1, either $L_0=L_1$
or $-k\ln(1+w)\leq L_1\leq 0\leq L_0$.
Thus, using (\ref{concave}), we obtain that with probability 1
\begin{align}
\nonumber
0\leq f(L_0)-f(L_1)
&\leq (L_0-L_1)
\frac{f(0)-f\big(-k\ln(1+w)\big)}{0-(-k\ln(1+w))}
\\
\nonumber
&=
(L_0-L_1)
\frac{w}{1+w}
\frac{-\ln(1+w)+\ln(1+\lambda)}{k\ln(1+w)},
\\
\label{frel}
&=
(L_0-L_1)
\frac{w}{k(1+w)}
\left(\frac{\ln(1+\lambda)}{\ln(1+w)}-1\right),
\end{align}
where we have used
\begin{align*}
f(0)&=-w\ln(1+w)/(1+w)
\\ 
\intertext{and} 
f(-k\ln(1+w))&=
-\frac{w}{1+w}\ln\left(1+we^{k\ln(1+w)}\right)
\\
&=-\frac{w}{1+w}\ln\left(1+w(1+w)^k\right)
\\
&=-\frac{w}{1+w}\ln(1+\lambda).
\end{align*}

Combining (\ref{expect}) and (\ref{frel}),
we get
$0\leq \EE (L_0-L_1)\leq \rho\big[\EE (L_0-L_1)\big]$,
where
\[
\rho=\frac{w}{1+w}\left(\frac{\ln(1+\lambda)}{\ln(1+w)}-1\right).
\]
To obtain that $\EE (L_0-L_1)=0$, 
and hence that reconstruction is impossible,
it's enough that $\rho<1$.
But
\begin{align*}
\rho
&<\frac{w}{1+w}\frac{\ln(1+\lambda)}{\ln(1+w)}
\\
&\leq \ln(1+\lambda)\sup_{w>0}\frac{w}{(1+w)\ln(1+w)}
\\
&=\ln(1+\lambda).
\end{align*}
So certainly if $\lambda\leq e-1$, then $\rho<1$ as desired.
Also $\rho$ is continuous as a function of 
$\lambda$ (or $w$),
so 
the threshold value $\lambda_2(k)$ is in fact strictly greater than $e-1$,
and the proof of Theorem \ref{hardtheorem} is done.

\section*{Acknowledgements}
\label{sec:ack}
Thanks to Yuri Suhov for much encouragement 
throughout the course of this work,
and to Elchanan Mossel for valuable conversations during 
the Isaac Newton Institute programme on Computation,
Combinatorics and Probability in 2002.

\bibliographystyle{alpha}
\bibliography{jbm}
\label{sec:biblio}

\medskip

\noindent
\textsc{LIAFA,\\
CNRS et Universit\'e Paris 7,\\
case 7014,\\
2 place Jussieu,\\
75251 Paris Cedex 05\\
FRANCE}\\
\texttt{James.Martin@liafa.jussieu.fr}

\end{document}